\documentclass[reqno,12pt]{amsart}

\newtheorem{theorem}{Theorem}[section]
\newtheorem{lemma}[theorem]{Lemma}

\newtheorem{proposition}[theorem]{Proposition}

\theoremstyle{definition}

\theoremstyle{remark}
\newtheorem{remark}[theorem]{Remark}

\newcommand{\mysection}[1]{\section{#1}
\setcounter{equation}{0}}

\newcommand{\bR}{\mathbb R}

\newcommand{\bT}{\mathbb T}
\newcommand{\bZ}{\mathbb Z}

\renewcommand{\epsilon}{\varepsilon}

\begin{document}
\title[quasi-geostrophic equations
]{Dissipative
quasi-geostrophic equations in critical Sobolev spaces: smoothing effect and global well-posedness}

\author[H. Dong]{Hongjie Dong}
\address[H. Dong]
{Division of Applied Mathematics, Brown University,
182 George Street, Providence, RI 02912, USA}
\email{hjdong@brown.edu}
\thanks{The author was partially supported by a
start-up funding from the Division of Applied Mathematics of Brown
University and the National Science Foundation under agreement No. DMS-0111298 and DMS-0800129.}



\begin{abstract}
We study the critical and super-critical  dissipative
quasi-geostrophic  equations in $\bR^2$ or $\bT^2$. An optimal local smoothing effect of solutions with arbitrary initial data in $H^{2-\gamma}$ is proved. As a main application, we establish the
global well-posedness for the critical 2D quasi-geostrophic
equations with periodic $H^1$ data. Some decay in time estimates are
also provided.
\end{abstract}

\maketitle

\mysection{Introduction}
                                                    \label{intro}

We are interested in the initial value problem of two dimensional
dissipative quasi-geostrophic equations
\begin{equation}
                                        \label{qgeq1}
\left\{\begin{array}{l l}
\theta_t+u\cdot \nabla\theta+(-\Delta)^{\gamma/2}\theta=0\quad &\text{on}\, \bR^2\times (0,\infty),\\
\theta(0,x)=\theta_0(x)\quad & x\in \bR^2,\end{array}\right.
\end{equation}
where $\gamma\in (0,2]$ is a fixed parameter and the velocity
$u=(u_1,u_2)$ is divergence free and is determined by the Riesz
transforms of the potential temperature $\theta$:
$$
u=(-{\mathcal R}_2\theta,{\mathcal R}_1\theta)=(-\partial_{x_2}(-\Delta)^{-1/2}\theta, \partial_{x_1}(-\Delta)^{-1/2}\theta).
$$

Equation \eqref{qgeq1} is an important model in geophysical fluid
dynamics. It is derived from general quasi-geostrophic equations in
the special case of constant potential vorticity and buoyancy
frequency. 
Mathematically, the equation has also been considered to be a 2D model of the 3D incompressible Navier-Stokes equations. It is therefore an
interesting model for investigating existence issues on genuine 3D
Navier-Stokes equations. Recently, this equation has been studied by
many authors, see \cite{cordoba,const1,const2,ju1,ju2,wu1,wu2,wu3} and
references therein.

The global existence of a weak solution to \eqref{qgeq1} follows from
Resnick \cite{resnick}. The cases $\gamma>1, \gamma=1$ and
$\gamma<1$ are  called sub-critical, critical and super-critical
respectively. The sub-critical case is well understood. Wu
established in \cite{wu1} the global existence of a unique regular
solution to \eqref{qgeq1} with initial data $\theta_0$ in $L^p$ for
$p>2/(\gamma-1)$. With initial data in the scaling invariant space
$L^{2/(\gamma-1)}$, the proof of the global well-posedness can be
found, for example, in recent \cite{Carrillo}, where the asymptotic
behavior of the solutions is also studied. By using a Fourier
splitting method, Constantin and Wu \cite{const2} showed the global
existence of a regular solution on the torus with periodic boundary
conditions and also a sharp $L^2$ decay estimate for weak solutions
with data in $L^2(\bR^2)\cap L^1(\bR^2)$. Furthermore, very recently
in \cite{dongli} the author and Li estimated the higher order
derivatives of the solution and proved that it is actually spatial
analytic.

However, the cases of critical and super-critical quasi-geostrophic
equations still have  quite a few unsolved problems. In the critical
case, Constantin, C\'ordoba and Wu \cite{const1} gave a construction
of global regular solutions for the initial data in $H^1$ under a
smallness assumption of $L^\infty$ norm of the data. Moreover, they
showed that the solutions are spatial analytic for sufficiently
large $t$. In Chae and Lee \cite{chae}, the global
existence and uniqueness were obtained for small initial data in the
critical Besov space $B_{2,1}^{2-\gamma}$.  In \cite{ju1}, Ju
improved Chae and Lee's result by showing that \eqref{qgeq1} is
globally well-posed for small data in $H^s$ if $s\geq 2-\gamma$, and
locally well-posed for large data if $s>2-\gamma$.

Very recently, there are two important papers \cite{kiselev} and \cite{caf}.
In \cite{kiselev} the global well-posedness for the critical quasi-geostrophic
equations with periodic $C^\infty$ data was established by Kiselev,
Nazarov and Volberg by proving certain non-local
maximum principle. In \cite{caf} Caffarelli and Vasseur constructed a global regular solution for the critical quasi-geostrophic equations
with $L^2$ initial data. To the best of our knowledge, the uniqueness of such weak solution is still open.

For results with minimal regularity assumptions, in recent \cite{miura}, Miura improved the result in \cite{ju1} and
proved the local in time existence of a unique regular solution for
large initial data in the critical Sobolev space $H^{2-\gamma}$. A similar result was also obtained independently in
Ju \cite{ju3} by using a different approach. For other results about
the critical and super-critical cases, we also refer the readers to
\cite{chae,CW06-Hol,ju1,ju2,wu2,wu3,yu07}.

Next we shall describe the main results of the present paper.

Our
first result (Theorem \ref{thm1} and \ref{thm3}) is concerning the {\em optimal} local smoothing effect of solutions. It says, roughly
speaking, that the smoothing effect of the equations in spaces  is
the same for the corresponding linear equations.  We remark that in the critical or super-critical cases, one
has higher derivative in the flow term $u\cdot\nabla \theta$ than in the
dissipation term $(-\Delta)^{\gamma/2}\theta$. A general
understanding is that the former term tends to make the smoothness
of $\theta$ worse, while the latter term tends to make it better. We
show that for small $\gamma\in (0,1]$, the dissipation is still
strong enough to balance the nonlinear term.  This result implies,
in particular, that the solution is infinitely differentiable. For
the critical quasi-geostrophic equation, although we have $H^1$
local well-posedness, to get global existence the authors of
\cite{kiselev} have to assume that the initial data is smooth. In
this connection, we note that the chief purpose of the current
article is to fill in this gap.

As a main application of Theorem
\ref{thm1}, in the second result (Theorem \ref{thmglobal}), we
obtain the global well-posedness of the critical quasi-geostrophic
equation with period $H^1$ data. We remark that the problem of $H^1$
global well-posedness of the critical quasi-geostrophic equation has
been open for years.\footnote{After the paper was finished, the
author and Dapeng Du realized that by adapting a method in
\cite{kiselev} with suitable modifications, the results here can be
used to establish the $H^1$ global well-posedness of the critical
quasi-geostrophic equation in the whole space. We present this in a
subsequent paper \cite{dongdu2}.} Moreover, we prove an exponential
decay estimate of the solution and all its derivatives, and show
that the solution is spatial analytic for large $t$. Although some
results here are based on the main result of \cite{miura}, the proof
of which uses the contraction argument, this article is
not a simple extension of \cite{miura}. The contraction argument
is not sufficient to establish the infinite differentiability
of the solution, since the time of existence of the solution in
$H^\beta,\beta\ge 1$
may be dependent on $\beta$. Instead, a suitable arrangement of the
nonlinear term enable us to use a bootstrap argument to get
the infinite differentiability as well as an exponential decay estimate.

In a forthcoming article, we are going to generalize these results to
more general Besov spaces. Although the main idea is similar, more
complicated arguments and estimates are involved.

The remaining part of the article is organized as follows: our main
theorems (Theorem \ref{thm1}, \ref{thm2}, \ref{thm3} and
\ref{thmglobal}) are stated in the next section. We define some
notation which we shall use later and recall some basic estimates in
Section \ref{linear}. The proof of a commutator estimate (Lemma
\ref{lem7.00}) is deferred to Section \ref{comm}. These estimate
enable us to prove Theorem \ref{thm1} and \ref{thm2} in Section
\ref{proofofthm1} by adapting an idea, which has been used in
\cite{dongdu,dongli,sawada,natasa}.
Section \ref{proofofthm3} and \ref{proofofglobal} are devoted to the
proofs of a Theorem \ref{thm3} and \ref{thmglobal}.

\mysection{Main theorems}
                                    \label{main}
Define $G(t,x)=G_\gamma(t,x)$ by its Fourier  transform
$\widehat{G_\gamma}(t,\xi)=e^{-t|\xi|^\gamma}$ for $t>0$. Then
$G_\gamma(t,x)$ is the fundamental solution of the linear operator
$\partial_t+(-\Delta)^{\gamma/2}$. It also has the scaling property
$$
G_\gamma(t,x)=t^{-\frac{2}{\gamma}}G_\gamma(1,xt^{-\frac{1}{\gamma}}).
$$
It is well-known that \eqref{qgeq1} can be rewritten into an
integral equation
$$
\theta(t,\cdot)=G(t,\cdot)*\theta_0-\int_0^t G(t-s,\cdot)*(u\cdot\nabla \theta)(s,\cdot)\,ds.
$$
Since $u$ is divergence free, integration by parts yields
$$
\theta(t,\cdot)=G(t,\cdot)*\theta_0-\int_0^t \nabla G(t-s,\cdot)*(u\theta)(s,\cdot)\,ds.
$$
In the sub-critical case, after obtaining  suitable linear and
bilinear estimates in certain Banach spaces, one can use the
classical Kato's contraction method \cite{fuKato} to prove the local
existence results. However, due to the weak dissipations, this
method seems not applicable in the usual way for the critical and
super-critical cases. In particular, it is difficult to find a
suitable Banach space $X$ so that the bilinear term is continuous
from $X\times X$ to $X$.

The following theorem is recently proved in Miura  \cite{miura} by
using a variation of the Kato's method combined with a commutator
estimate associated with the Littlewood-Paley operator in the
Sobolev space (see also recent Ju \cite{ju3} for a different
approach).
\begin{proposition}
            \label{prop0}
Let $\gamma\in (0,1]$ and $\theta_0\in H^{2-\gamma}$. Then
there exists $T>0$ such that the initial value problem for
\eqref{qgeq1} has a unique solution
\begin{equation}
                                            \label{eq4.18}
\theta(t,x)\in C([0,T]; H^{2-\gamma})\cap L^2(0,T; H^{2-\gamma/2}).
\end{equation}
The solution $\theta$ satisfies
\begin{equation}
            \label{eq2.08}
\sup_{0<t<T}t^{\beta/\gamma}\|\theta(t,\cdot)\|_{\dot H^{2-\gamma+\beta}}<\infty,
\end{equation}
for any $\beta\in [0,\gamma)$ and
\begin{equation}
            \label{eq2.10}
\lim_{t\to 0}t^{\beta/\gamma}\|\theta(t,\cdot)\|_{\dot H^{2-\gamma+\beta}}=0,
\end{equation}
for any $\beta\in (0,\gamma)$. Furthermore, there  exists
$\epsilon_0>0$ such that if $\|\theta_0\|_{\dot
H^{2-\gamma}}<\epsilon_0$, then we can take $T=\infty$.
\end{proposition}


By adapting the idea which were used in \cite{dongdu,sawada,dongli,natasa,MS2006}, we are
able to get the optimal local smoothing effect of the solution. Next we state our
main results.
\begin{theorem}
                    \label{thm1}
Let $\gamma\in (0,1]$ and $\theta_0\in H^{2-\gamma}$. Then the
solution $\theta$ in Proposition \ref{prop0} satisfies
\begin{equation}
            \label{eq2.08b}
\sup_{0<t<T}t^{\beta/\gamma}\|\theta(t,\cdot)\|_{\dot H^{2-\gamma+\beta}}<\infty,
\end{equation}
for any $\beta\geq 0$ and
\begin{equation}
            \label{eq2.10b}
\lim_{t\to 0}t^{\beta/\gamma}\|\theta(t,\cdot)\|_{\dot H^{2-\gamma+\beta}}=0,
\end{equation}
for any $\beta>0$.
\end{theorem}

\begin{remark}
                                    \label{rem3.52}
If we assume $\theta_0\in H^{2-\gamma}$, the Sobolev embedding
theorem, the boundedness of Riesz transforms on $L^p,1<p<\infty$
and Theorem \ref{thm1}
together with the $L^p$ maximum principle imply that the solution
$\theta$ and $u$ are smooth in $x$ in $(0,T)\times \bR^2$. Then from
the equation \eqref{qgeq1} itself, we see that they are also smooth
in $t$ in that region. Consequently, the mild solution $\theta$ is
in fact a classical solution of \eqref{qgeq1}.
\end{remark}

The proof of Theorem \ref{thm1} also yields an optimal decay in  time
estimate of higher order Sobolev norms in case of small initial
data.

\begin{theorem}
                    \label{thm2}
Under the assumptions of Theorem \ref{thm1}, there  exists
$\epsilon_0>0$ such that if $\|\theta_0\|_{\dot
H^{2-\gamma}}<\epsilon_0$, then

i) the initial value problem for \eqref{qgeq1} has a unique global
regular solution $\theta(t,x)$ in
$$C_b([0,\infty); H^{2-\gamma})\cap L^2((0,\infty);
H^{2-\gamma/2}).$$

ii) for any $\beta\geq 0$, the solution $\theta$ satisfies
\begin{equation}
            \label{eq2.08c}
\sup_{t>0}t^{\beta/\gamma}\|\theta(t,\cdot)\|_{\dot
H^{2-\gamma+\beta}}<\infty.
\end{equation}
\end{theorem}


Without much more work, a modification of the proof of Theorem
\ref{thm1} gives the integrability of the solution, along with its
derivatives, in time variable (See, e.g. \cite{dongdu}).
\begin{theorem}
                    \label{thm3}
Let $\gamma\in (0,1]$ and $\theta_0\in H^{2-\gamma}$. Then the
solution $\theta$ in Proposition \ref{prop0} satisfies
\begin{equation}
            \label{eq11.34}
\big\|t^{\beta_1/\gamma}\|\theta(t,\cdot)\|_{\dot
H^{2-\gamma+\beta}}\big\|_{L_t^{\gamma/\beta_2}(0,T)}<\infty,
\end{equation}
for any $\beta=\beta_1+\beta_2$ with $\beta_1\geq 0$ and $\beta_2\in
[0,\gamma/2]$.
\end{theorem}

\begin{remark}
As in Theorem \ref{thm2}, from the proof below  we can clearly see
that if the $\dot H^{2-\gamma}$ norm of the initial data is
sufficiently small (but independent of $\beta_1$ or $\beta_2$), then
one may take $T=\infty$ in Theorem \ref{thm3}.
\end{remark}

We can also consider the 2D quasi-geostrophic equations on the torus
with periodic boundary condition:
\begin{equation}
                                        \label{qgeq2}
\left\{\begin{array}{l l}
\theta_t+u\cdot\nabla\theta+(-\Delta)^{\gamma/2}\theta=0\quad
&\text{on}\, \bT^2\times (0,\infty),\\
\theta(0,x)=\theta_0(x)\quad & x\in \bT^2,\end{array}\right.
\end{equation}
where $\bT^2=[0,1]^2$ and $\theta_0\in \dot H^{2-\gamma}(\bT^2)$. As
usual, the zero-average condition is assumed:
$$
\int_{\bT^2}\theta_0(x)\,dx=0.
$$ Then by the Poincar\'e inequality, we have
$\theta_0\in H^{2-\gamma}(\bT^2)$. The proofs of Proposition
\ref{prop0} and Theorem \ref{thm1}, \ref{thm2}, \ref{thm3} can be
easily modified to get the corresponding results for \eqref{qgeq2}.
Also owing to a well-known fact $\int_{\bT^2}\theta(t,\cdot)\,dx=0$
and Poincar\'e's inequality, the homogeneous Sobolev norms in these
estimates can be replaced by the corresponding inhomogeneous norms.
We leave the details to interested readers.

For the critical quasi-geostrophic equations on the torus, we have
the following global existence result and {\em exponential} decay
estimate.
\begin{theorem}
                                        \label{thmglobal}
Let $\gamma=1$ and $\theta_0\in \dot H^{1}(\bT^2)$. Then  the
initial value problem for \eqref{qgeq2} has a unique global smooth
solution $\theta$ in
\begin{equation}
                                                    \label{eq22.11.21}
C_b([0,\infty); H^{1}(\bT^2))\cap L^2((0,\infty); H^{3/2}(\bT^2)).
\end{equation}
For some $T_0>0$, $\theta(t,\cdot)$ is spatial analytic for any
$t\geq T_0$. Furthermore, the solution and all its derivatives decay
exponentially as $t$ goes to infinity. More precisely, we have
\begin{equation}
            \label{eq22.10.59}
\sup_{t>0}e^{t/4}t^{\beta}\|\theta(t,\cdot)\|_{
H^{2-\gamma+\beta}}<\infty,
\end{equation} for any $\beta\geq 0$.
\end{theorem}

\mysection{Notation and some preliminary estimates}
                                            \label{linear}
First we recall the Littlewood-Paley decomposition. For any integer
$j$, define $\Delta_j$ to be the Littlewood-Paley projection
operator with $\Delta_j v=\phi_j*v$, where
$$
\hat\phi_j(\xi)=\hat\phi(2^{-j}\xi),\quad \hat\phi \in C_0^\infty(\bR^2\setminus\{0\}),\quad \hat\phi\geq 0,
$$
$$
\text{supp}\hat \phi\subset\{\xi\in \bR^2\,|\,1/2\leq |\xi|\leq 2\},
\quad \sum_{j\in \bZ}\hat\phi_j(\xi)=1\,\,\text{for}\,\,\xi\neq 0.
$$
Modulo a polynomials, formally we have the Littlewood-Paley decomposition
$$
v(\cdot,t)=\sum_{j\in \mathbb{Z}}\Delta_j v(\cdot,t).
$$
For any $p\in (1,\infty)$ and $s\geq 0$, as usual we denote $\dot
W^{s,p}$ and $W^{s,p}$ to be the homogeneous and inhomogeneous
Sobolev spaces with norms
$$
\|v\|_{\dot W^{s,p}}:=\Big\|(\sum_{k\in \mathbb{Z}}|2^{ks}\Delta_k
v|^2)^{1/2}\Big\|_{L^p}\sim \|\Lambda^s v\|_{L^p},
$$
$$
\|v\|_{W^{s,p}}:= \|v\|_{\dot W^{s,p}}+\|v\|_{L^p},
$$
with implicit constants depending on $p$ and $s$. When $p=2$, we use
$\dot H^s$ and $H^s$ instead of $\dot W^{s,p}$ and $W^{s,p}$.

Denote $\Lambda=(-\Delta)^{1/2}$. The following Bernstein's
inequality is well-known.

\begin{lemma}
                                        \label{bern}
Let $p\in [1,\infty]$ and $s\in \bR$. Then for any $j\in \bZ$, we
have
\begin{equation}
                                        \label{eq5.43}
\lambda 2^{js} \|\Delta_j v\|_{L^p}\leq \|\Lambda^s\Delta_j
v\|_{L^p}\leq \lambda' 2^{js} \|\Delta_j v\|_{L^p}
\end{equation}
with  some constants $\lambda$ and $\lambda'$ depending only on $p$
and $s$. Moreover, for $1\leq p\leq q\leq \infty$, there exists a
positive constant $C$ such that
\begin{equation}
                                        \label{eq5.48}
\|\Delta_j v\|_{L^q}\leq C2^{(2/p-2/q)j}\|\Delta_j v\|_{L^p}.
\end{equation}
\end{lemma}

We shall use the next two standard linear estimates, the proofs of
which can be found, for example, in \cite{miura}.
\begin{lemma}
                                        \label{lem2}
For any $\gamma>0$ and any function $v\in L^2$, we have
\begin{equation}
                                        \label{eq5.56}
e^{-2^{\gamma j+1}\lambda' t}\|\Delta_j v\|_{L^2}\leq
\|G(t,\cdot)*\Delta_j v\|_{L^2}\leq e^{-2^{\gamma j+1}\lambda
t}\|\Delta_j v\|_{L^2},
\end{equation}
where $\lambda$ and $\lambda'$ are some positive constants depending
only on $\gamma$.
\end{lemma}

\begin{lemma}
                                        \label{lem2.3}
For any $\gamma>0$ and $s\geq 0$, there exists a positive constant
$C$ depending only on $s$ and $\gamma$ such that for any $v\in L^2$,
we have
\begin{equation}
                                        \label{eq6.03}
\sup_{t\in (0,\infty)}t^{s/\gamma}\|G(t,\cdot)*v\|_{\dot H_x^s}\leq
C\|v\|_{L^2},
\end{equation}
\begin{equation}
                                        \label{eq6.05}
\lim_{t\to 0}t^{s/\gamma}\|G(t,\cdot)*v\|_{\dot H_x^s}=0,\quad \forall s>0.
\end{equation}
Moreover, for $s\in [0,\gamma/2]$ we have
\begin{equation}
                                        \label{eq6.06}
\|G(t,\cdot)*v\|_{L_t^{\gamma/s}\dot H_x^s}\leq C\|v\|_{L^2}
\end{equation}
\end{lemma}

As an easy consequence of Lemma \ref{lem2.3}, we have:
\begin{lemma}
                                        \label{lem2.4}
For any $\beta=\beta_1+\beta_2$ with $\beta_1\geq 0$ and $\beta_2\in
[0,\gamma/2]$, it holds that
\end{lemma}
\begin{equation}
            \label{eq6.14}
\big\|t^{\beta_1/\gamma}\|G(t,\cdot)*v\|_{\dot
H^{\beta}}\big\|_{L_t^{\gamma/\beta_2}}<C\|v\|_{L^2},
\end{equation}
where $C$ is a positive constant depending only on $\beta_1,\beta_2$
and $\gamma$.
\begin{proof}
By the semi-group property of the kernel $G(t,\cdot)$,
\eqref{eq6.03} and \eqref{eq6.06}, we get
\begin{align*}
&\big\|t^{\beta_1/\gamma}\|G(t,\cdot)*v\|_{\dot
H^{\beta}}\big\|_{L_t^{\gamma/\beta_2}}\\
&\quad =\big\|t^{\beta_1/\gamma}
\|G(t/2,\cdot)*G(t/2,\cdot)*v\|_{\dot
H^{\beta}}\big\|_{L_t^{\gamma/\beta_2}}\\
&\quad \leq C \|G(t/2,\cdot)*v\|_{L_t^{\gamma/\beta_2}\dot
H_x^{\beta_2}}\leq C\|v\|_{L^2}.
\end{align*}
The lemma is proved.
\end{proof}

The next lemma is a commutator estimate, which is a key estimate in
our proof. The proof of the lemma essentially follows that of
Proposition 2 \cite{miura}. We defer it to Section \ref{comm}.

\begin{lemma}
                                                    \label{lem7.00}
Assume $m\geq 0$, $1\leq s<2$, $t<1$ satisfying $m+t+s>0$. Then
there exists positive constant $C=C(s,t)$ such that
$$
\|[f,\Delta_j]g\|_{\dot H^m}\leq C2^{-(s+t-1)j}c_j(\|f\|_{\dot
H^{m+s}}\|g\|_{\dot H^{t}}+\|f\|_{\dot H^{s}}\|g\|_{\dot H^{m+t}})
$$ for any $j\in \bZ$, $f\in H^{m+t}$ and
$g\in H^{m+s}$ with $\|c_j\|_{l^2}\leq 1$.
Here,
$$[f,\Delta_j]g=f\Delta_j g-\Delta_j(fg).$$
\end{lemma}
\begin{remark}
                                             \label{rem7.05}
Define $\tilde \Delta_j=\sum_{|k-j|\leq 1}\Delta_j$. It is clear
from the proofs later that we only need a weaker estimate
\begin{align}
&\|\tilde \Delta_j[f,\Delta_j]g\|_{\dot H^m}\nonumber\\
                \label{eq5.38}
&\quad\quad \leq C2^{-(s+t-1)j}c_j(\|f\|_{\dot H^{m+s}}\|g\|_{\dot
H^{t}}+\|f\|_{\dot H^{s}}\|g\|_{\dot H^{m+t}}).
\end{align}
To get this estimate, the condition in Lemma \ref{lem7.00} can be
relaxed to $s<2$, $t<1$ and $m+t+s>0$.
\end{remark}

Finally, we shall also make use of the following lemma, which
follows simply from Plancherel's equality and localization property
of Littlewood-Paley projections. However, it is important in our
proofs.
\begin{lemma}
                                                    \label{lem8.51}
For any $j\in \bZ$ and $u,v\in L^2$, we have
\begin{equation}
                                                \label{eq8.53}
\int_{\bR^2} u\Delta_j v\,dx=\int_{\bR^2} (\tilde \Delta_j
u)(\Delta_j v)\,dx.
\end{equation}
\end{lemma}

\mysection{Local smoothing effect I}
                                        \label{proofofthm1}
Firstly, we give a general remark on our proofs. Recall that
equation \eqref{qgeq1} can be rewritten as
\begin{equation}
                    \label{eq7.18}
\theta(t,\cdot)=G(t,\cdot)*\theta_0-\int_0^t G(t-s,\cdot)*(u\cdot\nabla
\theta)(s,\cdot)\,ds.
\end{equation}
For the linear part, the estimate follows straightforwardly from
Lemma \ref{lem2.3} and \ref{lem2.4}. As usual, it is more difficult
to get a good estimate of the nonlinear term, especially in the
critical and super-critical case. Notice that the kernel
$G(t-s,\cdot)$ becomes singular as $s\to t$, and the initial data
$\theta_0$ is rough and only in $H^{2-\gamma}$. To deal with
the nonlinear term, the idea is to divide the integral into two
parts. For small $s$, we use the smoothness of the kernel
$G(t-s,\cdot)$. For large $s$ we should make use of the smoothness
of $\theta(s)$ and $u(s)$. This technique has been used in
\cite{dongdu,sawada}, and extensively in recent
\cite{dongli,natasa,MS2006}. Although the formulation
\eqref{eq7.18} does not appear explicitly in the proof below, we are
still able to exploit this idea. Moreover, thanks to the flexibility
of Lemma \ref{lem7.00}, the proof of the local smoothing effect is
considerably simpler comparing to those in \cite{dongdu,sawada,dongli,natasa,MS2006}.

However, since the estimates such as Bernstein's inequality and
fractional Leibniz's rule are quite rough, at present we are not
able to get any analyticity rate estimate as in \cite{sawada,dongli,MS2006}. On the other hand, it would be very
interesting to find out whether the mild solution of the critical
quasi-geostrophic equation with arbitrary $H^1$ initial data is spatially
analytic.
We note here that in the super-critical case even the solutions to
the corresponding linear equations are not spatially analytic. So
one should not expect that for the nonlinear equations.

{\bf Proof of Theorem \ref{thm1}:} Let $\theta$ be the solution in
Proposition \ref{prop0}. Denote $\theta_j=\Delta_j\theta$ and recall
$\Lambda=(-\Delta)^{1/2}$. For each $j\in \bZ$, we apply the
operator $\Delta_j$ to the both sides of \eqref{qgeq1} and get
$$
\partial_t \theta_j+\Delta_j(u\cdot \nabla
\theta)+\Lambda^\gamma \theta_j=0.
$$ Thus,
\begin{equation}
                                            \label{eq7.36}
\partial_t \theta_j+u\cdot \nabla
\theta_j+\Lambda^\gamma \theta_j=[u,\Delta_j]\nabla\theta.
\end{equation}
After multiplying both sides of \eqref{eq7.36} by $\theta_j$,
integrating in $x$ and noticing that $u$ is divergence free, we
obtain by using \eqref{eq5.43}, Lemma \ref{lem8.51} and H\"older's
inequality that
\begin{align*}
\frac{1}{2}\frac{d}{dt}\|\theta_j\|_{L^2}^2+\lambda 2^{\gamma j}
\|\theta_j\|_{L^2}^2 &\leq
\int_{\bR^2}\big([u,\Delta_j]\nabla\theta\big)\theta_j\,dx\\
&=\int_{\bR^2}\big(\tilde
\Delta_j[u,\Delta_j]\nabla\theta\big)\theta_j\,dx \\
&\leq
\big\|\tilde
\Delta_j[u,\Delta_j]\nabla\theta\big\|_{L^2}\|\theta_j\|_{L^2}.
\end{align*}
Therefore,
\begin{equation}
                                                \label{eq7.55}
\frac{d}{dt}\|\theta_j\|_{L^2}+\lambda 2^{\gamma
j}\|\theta_j\|_{L^2}\leq 2\big\|\tilde
\Delta_j[u,\Delta_j]\nabla\theta\big\|_{L^2}.
\end{equation}
Gronwall's inequality together with \eqref{eq7.55} yields
\begin{align}
&\|\theta_j(t,\cdot)\|_{L^2}\nonumber \\
                \label{eq8.05}
&\quad \leq e^{-2^{\gamma j}\lambda
t}\|\theta_j(0)\|_{L^2}+2\int_0^te^{-2^{\gamma j}\lambda (t-s)}
\big\|\tilde
\Delta_j[u,\Delta_j]\nabla\theta(s,\cdot)\big\|_{L^2}\,ds.
\end{align}

We prove the theorem by an induction on $\beta$. Proposition
\eqref{prop0} gives \eqref{eq2.08b} and \eqref{eq2.10b} for
$\beta\in (0,\gamma)$. Now assume $\beta_0\geq \gamma$, and
\eqref{eq2.08b} and \eqref{eq2.10b} are true for any $\beta\in
(0,\beta_0-\gamma/6]$. Let's consider the case when $\beta=\beta_0$.
We multiply the both sides of \eqref{eq8.05} by
$2^{(2-\gamma+\beta_0)j}$, use \eqref{eq5.43}, and split the
integral in to two parts,
\begin{equation}
                                            \label{eq11.00}
\|\theta_j(t,\cdot)\|_{\dot H^{2-\gamma+\beta_0}}\leq
2^{(2-\gamma+\beta_0)j}e^{-2^{\gamma j}\lambda
t}\|\theta_j(0,\cdot)\|_{L^2}+I_1+I_2,
\end{equation}
where
\begin{align*}
I_1&=2\int_0^{t/2}2^{(2-\gamma+\beta_0)j} e^{-2^{\gamma j}\lambda
(t-s)} \big\|\tilde
\Delta_j[u,\Delta_j]\nabla\theta(s,\cdot)\big\|_{L^2}\,ds,\\
I_2&=2\int^t_{t/2}2^{(2-\gamma+\beta_0)j} e^{-2^{\gamma j}\lambda
(t-s)} \big\|\tilde
\Delta_j[u,\Delta_j]\nabla\theta(s,\cdot)\big\|_{L^2}\,ds.
\end{align*}
We estimate $I_1$ and $I_2$ differently. In $I_1$, we absorb (most
part of)  the factor $2^{(2-\gamma+\beta_0)j}$ to the 'kernel'
$e^{-2^{\gamma j}\lambda (t-s)}$. While in $I_2$, we absorb (most
part of) that factor to the commutator term $\big\|\tilde
\Delta_j[u,\Delta_j]\nabla\theta(s,\cdot)\big\|_{L^2} $ and use the
localization property of $\tilde \Delta_j$ in the frequency space.

\textit{Estimate of $I_1$:} In Lemma \ref{lem7.00} we take
 $m=0$, $s=2-\gamma$,
 $t=1-3\gamma/4$, $f=u$, $g=\nabla \theta$, and get
\begin{align*}
I_1\leq C &c_j\int_0^{t/2}2^{(3\gamma/4+\beta_0)j}e^{-2^{\gamma
j}\lambda (t-s)}\|u(s,\cdot)\|_{\dot
H^{2-\gamma}}\|\theta(s,\cdot)\|_{\dot H^{2-3\gamma/4}}\,ds\\
\leq
C &c_j\int_0^{t/2}(t-s)^{-3/4-\beta_0/\gamma}\|\theta(s,\cdot)\|_{\dot
H^{2-\gamma}}\|\theta(s,\cdot)\|_{\dot H^{2-3\gamma/4}}\,ds\\
\leq C &c_j t^{-\beta_0/\gamma}\sup_{s\in
(0,t)}\|\theta(s,\cdot)\|_{\dot H^{2-\gamma}}\sup_{s\in (0,t)}
\big(s^{1/4}\|\theta(s,\cdot)\|_{\dot H^{2-3\gamma/4}}\big)\\
&\quad \cdot \int_{0}^{t/2}(t-s)^{-3/4}s^{-1/4}\,ds\\
\leq C &c_j t^{-\beta_0/\gamma}\sup_{s\in
(0,t)}\|\theta(s,\cdot)\|_{\dot H^{2-\gamma}}\sup_{s\in (0,t)}
\big(s^{1/4}\|\theta(s,\cdot)\|_{\dot H^{2-3\gamma/4}}\big),
\end{align*}
where in the second inequality we use the boundedness of Riesz
transforms in $L^2$.

\textit{Estimate of $I_2$:} By the Bernstein's inequality, it holds
that
$$
2^{kj}\|\tilde \Delta_j[u,\Delta_j]\nabla\theta(s)\big\|_{L^2}\leq
C\|\tilde\Delta_j[u,\Delta_j]\nabla\theta(s)\big\|_{\dot H^k}.
$$
Recall that here we assume $\beta_0\geq \gamma$. In Lemma
\ref{lem7.00} we take $m=\beta_0-\gamma/2\geq 0$, $s=2-2\gamma/3$,
$t=1-2\gamma/3$, $f=u$ and $g=\nabla \theta$, and get
\begin{align*}
I_2\leq C&\int^t_{t/2}2^{(2-\gamma/2)j} e^{-2^{\gamma j}\lambda
(t-s)} \big\|\tilde
\Delta_j[u,\Delta_j]\nabla\theta(s,\cdot)\big\|_{\dot H^{\beta_0-\gamma/2}}\,ds\\
\leq C &c_j\int_{t/2}^t 2^{5\gamma j/6}e^{-2^{\gamma j}\lambda
(t-s)}\big(\|u(s,\cdot)\|_{\dot
H^{2+\beta_0-7\gamma/6}}\|\theta(s,\cdot)\|_{\dot H^{2-2\gamma/3}}\\
&\quad +\|u(s,\cdot)\|_{\dot
H^{2-2\gamma/3}}\|\theta(s,\cdot)\|_{\dot
H^{2+\beta_0-7\gamma/6}}\big)\,ds\\
\leq C &c_j\int_{t/2}^t (t-s)^{-5/6}\|\theta(s,\cdot)\|_{\dot
H^{2+\beta_0-7\gamma/6}}\|\theta(s,\cdot)\|_{\dot
H^{2-2\gamma/3}}\,ds\\
\leq C &c_j\sup_{s\in
(0,t)}\big(s^{\beta_0/\gamma-1/6}\|\theta(s,\cdot)\|_{\dot
H^{2+\beta_0-7\gamma/6}}\big) \sup_{s\in
(0,t)}\big(s^{1/3}\|\theta(s,\cdot)\|_{\dot H^{2-2\gamma/3}}\big)\\
&\quad \cdot t^{-\beta_0/\gamma}\int_{t/2}^t (t-s)^{-5/6}s^{-1/6}\,ds\\
\leq C &c_j t^{-\beta_0/\gamma}\sup_{s\in
(0,t)}\big(s^{\beta_0/\gamma-1/6}\|\theta(s,\cdot)\|_{\dot
H^{2+\beta_0-7\gamma/6}}\big)\\
&\quad \cdot \sup_{s\in (0,t)}\big(s^{1/3} \|\theta(s,\cdot)\|_{\dot
H^{2-2\gamma/3}}\big),
\end{align*}
where in the second inequality we again use the boundedness of Riesz
transforms.



Now we take the $l_2$ norm of both sides of \eqref{eq11.00} in
$j\in\{-N,-N+1,\cdots,N-1,N\}$ for some positive integer $N$ and
then multiply both sides by $t^{\beta_0/\gamma}$. Owing to
\eqref{eq5.43} and Lemma \ref{lem2}, \ref{lem2.3}, it holds that
\begin{align}
&t^{\beta_0/\gamma} \Big(\sum_{j=-N}^N(\|\theta_j\|^2_{\dot
H_x^{2-\gamma+\beta_0}}\Big)^{1/2}\nonumber\\
&\quad \leq C\sup_{s\in
(0,C_1t)}\big(s^{\beta_0/\gamma}\|G(s,\cdot)*\theta_0\|_{\dot
H^{2-\gamma+\beta_0}}\big)\nonumber\\
&\quad\, +C \sup_{s\in (0,t)}\|\theta(s,\cdot)\|_{\dot
H^{2-\gamma}}\sup_{s\in (0,t)} \big(s^{1/4}\|\theta(s,\cdot)\|_{\dot
H^{2-3\gamma/4}}\big)\nonumber\\
                                            \label{eq11.59}
&\quad\, +C  \sup_{s\in
(0,t)}\big(s^{\frac{\beta_0}{\gamma}-\frac{1}{6}}
\|\theta(s,\cdot)\|_{\dot H^{2+\beta_0-\frac{7\gamma}{6}}}\big)
\sup_{s\in (0,t)}\big(s^{\frac{1}{3}}\|\theta(s,\cdot)\|_{\dot
H^{2-\frac{2\gamma}{3}}}\big),
\end{align}
where $C$ and $C_1$ are positive constants independent of $t$. In
the above inequality, the first term on the right-hand side is
bounded with respect to $t$ and goes to zero as $t\to 0$ due to
Lemma \ref{lem2.3}. The second and the third term is bounded for
$t\in (0,T)$ and go to zero as $t\to 0$ by the inductive assumption.
Letting $N\to +\infty$ in \eqref{eq11.59} yields \eqref{eq2.08b} and
\eqref{eq2.10b} for $\beta=\beta_0$. Theorem \ref{thm1} is then
proved.

{\bf Proof of Theorem \ref{thm2}:} The proofs of the first part of
the theorem and the second part for $\beta \in [0,\gamma)$ can be
found in \cite{miura}. We only need to show the second part for
$\beta\geq \gamma$. However, this follows immediately from the
induction argument in the proof of Theorem \ref{thm1} and
\eqref{eq11.59}. This completes the proof of Theorem \ref{thm2}.

 \mysection{Local smoothing effect II}
                                        \label{proofofthm3}
This section is devoted to the proof of Theorem \ref{thm3}. First we consider the case when $\beta\in [0,\gamma/2]$. As
$$
\theta \in L^\infty([0,T]; H^{2-\gamma})\cap L^2(0,T;
H^{2-\gamma/2}),
$$
for $\beta_1=0$ and $\beta_2=\beta\in [0,\gamma/2]$, by H\"older's
inequality and the interpolation estimate, we obtain
\begin{equation}
                                                    \label{eq12.33}
\theta \in L^{\gamma/\beta}([0,T]; \dot H^{2-\gamma+\beta}).
\end{equation}
This together with \eqref{eq2.08b} concludes Theorem \ref{thm3} in
its full generality when $\beta\in [0,\gamma/2]$.

Next we assume $\beta_0>\gamma/2$ and proceed by an induction on
$\beta$. Suppose \eqref{eq11.34} has been proved for $\beta\in
[0,\beta_0-\gamma/6]$. Let's consider the case when $\beta=\beta_0$
and assume $\beta_0=\beta_1+\beta_2$ for some $\beta_1>0$ and
$\beta_2\in [0,\gamma/2]$. Note that the estimates of both $I_1$ and
$I_2$ still holds true if we only assume $\beta_0>\gamma/2$. Because
of Theorem \ref{thm1}, we already know that $\theta(t,\cdot)\in \dot
H^{2-\gamma+\beta}$ for any $t>0$. Taking the $l_2$ norm of both
sides of \eqref{eq11.00} in $j\in \bZ$ and then multiply both sides
by $t^{\beta_1/\gamma}$ instead of $t^{\beta_0/\gamma}$ in the
previous section, we obtain
\begin{equation}
                                                \label{eq1.51}
t^{\beta_1/\gamma}\|\theta(t,\cdot)\|_{\dot
H^{2-\gamma+\beta_0}}\leq
Ct^{\beta_0/\gamma}\|G(C_1t,\cdot)*\theta_0\|_{\dot
H^{2-\gamma+\beta_0}}+CI_3+CI_4,
\end{equation}
where
\begin{align*}
I_3&=\int_0^{t/2}t^{\beta_1/\gamma}(t-s)^{-3/4-\beta_0/\gamma}
\|u(s,\cdot)\|_{\dot H^{2-\gamma}}\|\theta(s,\cdot)\|_{\dot
H^{2-3\gamma/4}}\,ds,\\
I_4&=\int_{t/2}^t t^{\beta_1/\gamma}
(t-s)^{-5/6}\|\theta(s,\cdot)\|_{\dot
H^{2+\beta_0-7\gamma/6}}\|\theta(s,\cdot)\|_{\dot
H^{2-2\gamma/3}}\,ds.
\end{align*}
We then show that all the three terms on the right-hand side of
\eqref{eq1.51} are in $L^{\gamma/\beta_2}(0,T)$.

Due to Lemma \ref{lem2.4}, the first term is indeed in
$L^{\gamma/\beta_2}(0,\infty)$. For $I_3$, we compute
$$
I_3\leq
C\int_0^{t/2}(t-s)^{-3/4-\beta_2/\gamma}\|\theta(s,\cdot)\|_{\dot
H^{2-\gamma}} \|\theta(s,\cdot)\|_{\dot H^{2-3\gamma/4}}\,ds.
$$
Owing to Proposition \ref{prop0}, we have
$$
\|\theta(t,\cdot)\|_{\dot H^{2-\gamma}}\|\theta(t,\cdot)\|_{\dot
H^{2-3\gamma/4}}\in L^4((0,T)).
$$
This together with the fractional integration yields
$$
I_3\in
L^{\gamma/\beta_2}(0,T).
$$
Finally, $I_4$ is less than
$$
C\int_{t/2}^t
(t-s)^{-5/6}\big(s^{\beta_1/\gamma-1/6}\|\theta(s,\cdot)\|_{\dot
H^{2+\beta_0-7\gamma/6}}\big)\big(s^{1/6}\|\theta(s,\cdot)\|_{\dot
H^{2-2\gamma/3}}\big)\,ds
$$
By the inductive assumption, we have,
\begin{align*}
t^{1/6}\|\theta(t,\cdot)\|_{\dot H^{2-2\gamma/3}}&\in L^6((0,T)),\\
t^{\beta_1/\gamma-1/4}\|\theta(t,\cdot)\|_{\dot
H^{2+\beta_0-5\gamma/4}}&\in L^{\gamma/\beta_2}((0,T)).
\end{align*}
These estimate together with the fractional integration yield
$I_4\in L^{\gamma/\beta_2}(0,T)$. It follows that \eqref{eq11.34}
holds for $\beta=\beta_0$. The theorem is proved.

\mysection{Global well-posedness when $\gamma=1$}
                                            \label{proofofglobal}

As we discussed in Remark \ref{rem3.52}, the solution $\theta$ and
$u$ become smooth immediately for $t>0$. Fix a $t_1\in (0,T)$. Then
we can consider $\theta(t_1)$ as initial data and apply the result
of the global existence for smooth initial data in \cite{kiselev}.
The boundedness of $\theta$ and its derivatives follows from the
uniform bound
\begin{equation}
                                            \label{eq11.55.22}
\|\nabla\theta(t,\cdot)\|_{L^\infty}\leq
C\|\nabla\theta_0\|_{L^\infty}e^{e^{C\|\theta_0\|_{L^\infty}}}
\end{equation}
established in \cite{kiselev} and Theorem \ref{thm1}. The solution
is in $C([0,\infty); H^{1})\cap L^2_{\text{loc}}((0,\infty);
H^{3/2})$. The uniqueness then follows in a standard way from the
local uniqueness result (see, e.g. \cite{miura}). To see the
solution is also in \eqref{eq22.11.21}, it suffices to verify the
decay estimate \eqref{eq22.10.59}.

Denote $\hat \theta(t,j),j\in \bZ^2$ to be the Fourier coefficients of
$\hat \theta(t,\cdot)$. Recall that $\theta(t,0)\equiv 0$ for any $t\geq
0$. Since $\theta$ and $u$ are smooth, Theorem 4.1 of C\'ordoba and
C\'ordoba \cite{cordoba} yields the following lemma.
\begin{lemma}
                                                    \label{lem11.26}
Under the assumptions of Theorem \ref{thmglobal}, there exists a
positive constant $C$ depending only on $\theta_0$ so that
$$
\|\theta(t,\cdot)\|_{L^\infty}\leq C/(1+t)
$$ for any $t\geq t_1$.
\end{lemma}

Thus we can choose $t$ large so that
$\|\theta(t,\cdot)\|_{L^\infty}$ is as small as we want. This
together with a small data result due to Constantin, C\'ordoba and
Wu \cite{const1} implies the spatial analyticity of $\theta$ for
$t\geq T_0$ for some $T_0\geq t_1$. More precisely, we have
\begin{lemma}
                                                    \label{lem11.43}
Under the assumptions of Theorem \ref{thmglobal}, there exists
$T_0\geq t_1$ such that
\begin{equation}
                                                \label{eq11.45}
y(t):=\sum_{j\in \bZ^2\setminus\{(0,0)\}}|\hat
\theta(t,j)|e^{(t-T_0)|j|/2}\leq 1/2,
\end{equation} for any $t\geq T_0$.
\end{lemma}

We claim that \eqref{eq11.45} implies \eqref{eq22.10.59}. Indeed,
for $t\in (0,T_0)$ estimate \eqref{eq22.10.59} is an immediate
consequence of Theorem \ref{thm1}, \eqref{eq11.55.22} and
Poincar\'e's inequality. For any $t\geq T_0$, we have
\begin{align*}
e^{t/2}t^{2\beta/\gamma}\|\theta(t,\cdot)\|_{
H^{2-\gamma+\beta}}^2&\leq Ce^{t/2}t^{2\beta/\gamma}\sum_{j\in
\bZ^2\setminus\{(0,0)\}}|\hat
\theta(t,j)|^2|j|^{2(2-\gamma+\beta)}\\
&\leq C(y(t))^2\leq C.
\end{align*}
This finishes the proof of Theorem \ref{thmglobal}.

\mysection{A commutator estimate}
                                        \label{comm}
This section is devoted to the proof of Lemma \ref{lem7.00}. We
follow closely the idea of Proposition 2 in \cite{miura}  (see
also earlier \cite{chae,chemin,danchin} for
similar estimates). However, since we also consider higher order
Sobolev norms by introducing a parameter $m$, we give a proof here
for the sake of completeness. It is worth noting that from the proof
below the condition of Lemma \ref{lem7.00} can be relaxed.

We start with the definition of Bony's paraproduct operator and some
basic estimates for the paraproduct operator (see, e.g.
\cite{runst}). Define paraproduct operators by
$$
T_f g:=\sum_{j\in \bZ}S_j f\Delta_j g,\quad R(f,g):=\sum_{|i-j|\leq 2}\Delta_i f \Delta_j g,
$$
where $S_j f=\sum_{k\leq j-3}\Delta_k f.$ Then we have
$$
fg=T_f g+T_g f+R(f,g).
$$

\begin{lemma}
                    \label{lem2.01}
i) If $s<1,t\in \bR$, there exists a positive constant $C$ depending only on $s$ and $t$ such that for any $f\in \dot H^s(\bR^2)$ and $g\in \dot H^t(\bR^2)$ we have
\begin{equation}
            \label{eq2.04}
\|T_f g\|_{\dot H^{s+t-1}}\leq C\|f\|_{\dot H^s}\|g\|_{\dot H^t}.
\end{equation}

ii) If $s+t>0$, there exists a positive constant $C$ depending only on $s$ and $t$ such that for any $f\in \dot H^s(\bR^2)$ and $g\in \dot H^t(\bR^2)$ we have
\begin{equation}
            \label{eq2.07}
\|R(f,g)\|_{\dot H^{s+t-1}}\leq C\|f\|_{\dot H^s}\|g\|_{\dot H^t}.
\end{equation}

iii) If $s,t<1$ and $s+t>0$, there exists a positive constant $C$ depending only on $s$ and $t$ such that for any $f\in \dot H^s(\bR^2)$ and $g\in \dot H^t(\bR^2)$ we have
\begin{equation}
            \label{eq2.088}
\|fg\|_{\dot H^{s+t-1}}\leq C\|f\|_{\dot H^s}\|g\|_{\dot H^t}.
\end{equation}
\end{lemma}

The following fractional Leibniz's rule is well-known.

\begin{lemma}
                    \label{lem2.56}
Assume $s\geq 0$ and $p\in (1,\infty)$. Then we have
$$
\|fg\|_{\dot W^{s,p}}\leq C\|f\|_{\dot W^{s,p_1}}
\|g\|_{L^{p_2}}+C\|f\|_{L^{p_1'}}\|g\|_{\dot W^{s,p_2'}}
$$ if the right-hand side is finite. Here $p_1,p_2,p_1',p_2'\in (1,+\infty)$ satisfy
$$
1/p=1/p_1+1/p_2=1/p_1'+1/p_2'.
$$
\end{lemma}

Now we are ready to prove Lemma \ref{lem7.00}. Denote
$$
f_j=\Delta_j f,\quad g_j=\Delta_j g,\quad \check
\Delta_j=\sum_{|k-j|\leq 2}\Delta_k
$$ for any $j\in \bZ$. In terms of paraproducts, we have
\begin{align*}
[f,\Delta_j] g&=-\Delta_j R(f,g)-\Delta_j(T_g f)+[T_f,\Delta_j]g+R(f,g_j)+T_{g_j} f\\
&=-\Delta_j R(f,g)-\Delta_j(T_g f)+\sum_{|k-j|\leq 3}[S_k f,\Delta_j]g_k\\
&\quad +\sum_{|k-j|\leq 2}\check \Delta_k f\Delta_k g_j+\sum_{k\geq j+1}S_kg_j f_k\\
&:=I_1+I_2+I_3+I_4+I_5,
\end{align*}
where in the second  equality above we use the  localization
property of Littlewood-Paley projections in the frequency space.
Choose $p_1\in (2,\infty)$ sufficiently large so that $s+2/p_1<2$.
This is possible because $s<2$. Let $p_2\in (2,\infty)$ be a number
satisfying $1/p_1+1/p_2=1/2$.

\textit{Estimate of $I_1$:} Because $m+s+t>0$, by using \eqref{eq2.07} with $m+s$ and $t$ in place of $s$ and $t$ respectively, we get
$$
\|I_1\|_{\dot H^m}\leq Cc_j 2^{(1-s-t)j}\|f\|_{\dot H^{m+s}}\|g\|_{\dot H^t},
$$
where
$$
c_j=2^{(s+t-1)j}\|\Delta_j R(f,g)\|_{\dot H^m}/\| R(f,g)\|_{\dot H^{m+s+t-1}}.
$$

\textit{Estimate of $I_2$:} Since $t<1$, \eqref{eq2.04} with $t+m$ in place of $t$ gives
$$
\|I_2\|_{\dot H^m}\leq \tilde c_j 2^{(1-s-t)j}\|f\|_{\dot H^{s+m}}\|g\|_{\dot H^{t}},
$$
where
$$
\tilde c_j=2^{(s+t-1)j}\|\Delta_j T_g f\|_{\dot H^m}/\|T_g f\|_{\dot H^{m+s+t-1}}
$$

\textit{Estimate of $I_3$:} The estimate of $I_3$ is more delicate. By the mean value theorem, we have
\begin{align*}
I_3 &=\sum_{|k-j|\leq 3}\int_{\bR^2}\int_0^1 \phi_j(y)y(S_k\nabla f)(x-sy)g_k(x-y)\,dsdy\\
&=2^{-j}\sum_{|k-j|\leq 3}\int_{\bR^2}\int_0^1 \phi(y)y(S_k\nabla f)(x-2^{-j}sy)g_k(x-2^{-j}y)\,dsdy.
\end{align*}
Now due to Minkowski's inequality and Lemma \ref{lem2.56}, we get
\begin{equation}
                    \label{eq3.41}
\|I_3\|_{\dot H^m}\leq C2^{-j}\sum_{|k-j|\leq 3}\Big(\|S_k\nabla f\|_{\dot W^{m,p_1}}\|g_k\|_{L^{p_2}}+
\|S_k\nabla f\|_{L^{p_1}}\|g_k\|_{\dot W^{m,p_2}}\Big).
\end{equation}
Recall $s+2/p_1<2$. Then by H\"older's inequality,
$$
|S_k\Lambda^m\nabla f|\leq C  2^{(2-s-2/p_1)k}\|2^{(s+2/p_1-2)i}\Lambda^m\nabla f_i\|_{l^2}.
$$
Therefore,
\begin{align*}
\|S_k\nabla f\|_{\dot W^{m,p_1}}&\leq C  2^{(2-s-2/p_1)k}\big\|\|2^{(s+2/p_1-2)i}\Lambda^m\nabla f_i\|_{l^2}\big\|_{L^{p_1}}\\
&=C2^{(2-s-2/p_1)k}\|f\|_{\dot W^{m+s+2/p_1-1,p_1}}\\
&\leq C2^{(2-s-2/p_1)k}\|f\|_{\dot H^{m+s}},
\end{align*}
where in the last inequality we use Sobolev embedding theorem. Similarly,
$$
\|S_k\nabla f\|_{L^{p_1}}\leq C2^{(2-s-2/p_1)k}\|f\|_{\dot H^{s}}.
$$
These estimates together with \eqref{eq3.41} and Lemma \ref{bern} yield
\begin{align*}
&\|I_3\|_{\dot H^m}\\
&\quad \leq C 2^{(1-s-2/p_1)j}\Big(\|f\|_{\dot H^{m+s}}\sum_{|k-j|\leq 3}\|g_k\|_{L^{p_2}}+\|f\|_{\dot H^{s}}
\sum_{|k-j|\leq 3}\|g_k\|_{\dot W^{m,p_2}}\Big)\\
&\quad\leq C 2^{(1-s-t)j}\Big(\|f\|_{\dot H^{m+s}}\sum_{|k-j|\leq 3}2^{(t-2/p_1)k}\|g_k\|_{L^{p_2}}\\
&\quad\quad +\|f\|_{\dot H^{s}}
\sum_{|k-j|\leq 3}2^{(t-2/p_1)k}\|g_k\|_{\dot W^{m,p_2}}\Big)\\
&\quad\leq C 2^{(1-s-t)j}\Big(\|f\|_{\dot H^{m+s}}\sum_{|k-j|\leq 3}2^{tk}\|g_k\|_{L^{2}}+\|f\|_{\dot H^{s}}
\sum_{|k-j|\leq 3}2^{tk}\|g_k\|_{\dot H^{m}}\Big)\\
&\quad\leq C2^{(1-s-t)j}\bar c_j(\|f\|_{\dot H^{m+s}}\|g\|_{\dot H^{t}}+\|f\|_{\dot
H^{s}}\|g\|_{\dot H^{m+t}}),
\end{align*}
where
$$
\bar c_j=\frac{1}{100}\sum_{|k-j|\leq 3}2^{tk}\|g_k\|_{L^{2}}/\|g\|_{\dot H^t}+\frac{1}{100}\sum_{|k-j|\leq 3}2^{tk}\|g_k\|_{\dot H^{m}}/\|g\|_{\dot H^{m+t}}.
$$ It is easily seen that $\|\bar c_j\|_{l^2}\leq 1$, which completes the estimate of $I_3$.

\textit{Estimate of $I_4$:} Lemma \ref{lem2.56} yields
\begin{align*}
\|I_4\|_{\dot H^m} &\leq \sum_{|k-j|\leq 2}\|\check \Delta_k f\Delta_k g_j\|_{\dot H^m}\\
&\leq C\sum_{|k-j|\leq 2}\Big(\|\check \Delta_k f\|_{\dot W^{m,p_1}}\|\Delta_k g_j\|_{L^{p_2}}+\|\check \Delta_k f\|_{L^{p_1}}\|\Delta_k g_j\|_{\dot W^{m,p_2}}\Big)\\
&\leq C2^{(1-s-t)j}\sum_{|k-j|\leq 2}\Big(2^{k(s+2/p_1-1)}\|\check \Delta_k f\|_{\dot W^{m,p_1}}2^{k(t-2/p_1)}\|\Delta_k g_j\|_{L^{p_2}}\\
&\quad +2^{k(s+2/p_1-1)}\|\check \Delta_k f\|_{L^{p_1}}2^{k(t-2/p_1)}\|\Delta_k g_j\|_{\dot W^{m,p_2}}\Big)\\
&\leq C2^{(1-s-t)j}\bar c_j(\|f\|_{\dot H^{m+s}}\|g\|_{\dot H^{t}}+\|f\|_{\dot
H^{s}}\|g\|_{\dot H^{m+t}}),
\end{align*}
where $\bar c_j,j\in \bZ$ are the same constants  as in the estimate
of $I_3$.

\textit{Estimate of $I_5$:}
By using the boundedness of the operator $S_k$ in $L^p,p\in (1,\infty)$ and Lemma \ref{lem2.56} , we have
\begin{align*}
&\|I_5\|_{\dot H^m}\\
&\,\leq C\|g_j\|_{\dot W^{m,p_2}}\sum_{k\geq j+1}\|f_k\|_{L^{p_1}}
+C\|g_j\|_{L^{p_2}}\sum_{k\geq j+1}\|f_k\|_{\dot W^{m,p_1}}\\
&\,:= I_{51}+I_{52}.
\end{align*}
By Lemma \ref{bern} and H\"older's inequality,
\begin{align*}
I_{51}&\le C2^{(-t+1-2/p_2)j}\|g_j\|_{\dot H^{m+t}}\sum_{k\geq j+1}
2^{(1-2/p_1)k}\|f_k\|_{L^2}\\
&\leq C\bar c_j2^{(-t+1-2/p_2)j}\|g\|_{\dot H^{m+t}}
\left(\sum_{k\geq j+1}2^{2sk}\|f_k\|_{L^2}^2\right)^{1/2}2^{(1-2/p_1-s)j}\\
&\leq C\bar c_j2^{(-t-s+1)j}\|g\|_{\dot H^{m+t}}\|f\|_{\dot H^s}.
\end{align*}
In a similar way,
\begin{equation*}
I_{52}\leq C\bar c_j2^{(-t-s+1)j}\|g\|_{\dot H^{t}}\|f\|_{\dot H^{m+s}}.
\end{equation*}

Combining all these estimates together finishes the proof of the
lemma.

As we mentioned in Remark \ref{rem7.05}, in the proofs of the main
theorems we only use the estimate of a frequency localized object
$\tilde \Delta_j[f,\Delta_j] g$ instead of $[f,\Delta_j] g$ itself.
Notice that
$$
\tilde \Delta_j I_5=\tilde \Delta_j \sum_{j+1 \leq k\leq j+4}S_kg_j f_k.
$$ Now due to the finiteness of the  number of the sum on $k$
and boundedness of $\tilde \Delta_j$, in the estimate of $\tilde
\Delta_j I_5$ the condition that $s\geq 1$ can be removed. Moreover,
in the estimates of $I_3$, $I_4$ and $I_5$, where Lemma
\ref{lem2.56} is applied, we may estimate $ \|\tilde\Delta_j
I_l\|_{\dot H^{m+s+t}},l=3,4,5$ instead of $\|\Delta_j I_l\|_{\dot
H^{m}},l=3,4,5$ and still get the same bounds. Therefore, the
condition $m\geq 0$ can also be removed too. Since these are the
only places using these two conditions, we remark that to obtain
\eqref{eq5.38} we only require $s<2$, $t<1$ and $m+t+s>0$.

\section*{Acknowledgment}
The author would like to thank Peter Constantin, Zhen Lei, Dong Li and the referees for their very helpful comments.

\end{document}